\documentclass[reqno]{amsart}

\usepackage{amsmath,amssymb,amsthm}
\usepackage{caption}

\usepackage{graphicx,tikz}
\newtheorem{theorem}{Theorem}
\newtheorem*{thm}{Theorem}

\newtheorem{lemma}{Lemma}

\theoremstyle{definition}

\theoremstyle{remark}

\DeclareMathOperator{\rk}{rank}

\begin{document}

\title[]{Approximate Solutions of Linear Systems \\at a universal rate}
\thanks{S.S. is supported by the NSF (DMS-2123224) and the Alfred P. Sloan Foundation.}

\author[]{Stefan Steinerberger}
\address{Department of Mathematics, University of Washington, Seattle}
\email{steinerb@uw.edu}

\begin{abstract} Let $A \in \mathbb{R}^{n \times n}$ be invertible, $x \in \mathbb{R}^n$ unknown and $b =Ax $ given. We are interested in \textit{approximate} solutions: vectors $y \in \mathbb{R}^n$ such that $\|Ay - b\|$ is small. We prove that for all $0< \varepsilon <1 $ there is a composition of $k$ orthogonal projections onto the $n$ hyperplanes generated by the rows of $A$, where
$$k \leq    2 \log\left(\frac{1}{\varepsilon} \right) \frac{ n}{ \varepsilon^{2}}$$
which maps the origin to a vector $y\in \mathbb{R}^n$ satisfying $\| A y - Ax\| \leq \varepsilon \cdot \|A\| \cdot \| x\|$. We note that this upper bound on $k$ is independent of the matrix $A$.
This procedure is stable in the sense that $\|y\| \leq 2\|x\|$. The existence proof is based on a probabilistically refined analysis of the Random Kaczmarz method which seems to achieve this rate when solving for $A x = b$ with high likelihood. \end{abstract}

\maketitle

\vspace{-8pt}

\section{Introduction and Main Result}
\subsection{Setup.} Throughout this paper, let $A \in \mathbb{R}^{m \times n}$, $m \geq n$, be a (possibly overdetermined) linear system of equations, suppose that $A$ is injective and that
$ Ax =b,$
where $x \in \mathbb{R}^n$ is the (unknown) unique solution and $b$ is a given right-hand side.
We use $a_1, \dots, a_m \in \mathbb{R}^n$ to denote the rows of $A$, $\|A\|$ to denote the operator norm, $\|A\|_F$ for the Frobenius norm and $\sigma_{\min}$ for the smallest singular value.
The linear system $Ax = b$ can be written as
$$ \left\langle a_i, x \right\rangle = b_i \qquad \mbox{for}~1 \leq i \leq m,$$
one can interpret $x$ as the unique point at the intersection of the hyperplanes 
$$ H_i = \left\{w \in \mathbb{R}^n: \left\langle w, a_i \right\rangle = b_i \right\}  \qquad \mbox{for}~1 \leq i \leq m.$$

\begin{center}
\begin{figure}[h!]
\begin{tikzpicture}[scale=2.5]
\draw[thick] (-0.5, -0.5) -- (1, 1);
\draw[thick] (0, -0.5) -- (0, 1);
\draw[thick] (0.5, -0.5) -- (-1, 1);
\filldraw (0,0) circle (0.03cm);
\node at (0.2, 0) {$x$};
\filldraw (0.7, 0.7) circle (0.03cm);
\node at (0.9, 0.7) {$y$};
\filldraw (0, 0.7) circle (0.03cm);
\node at (0.25, 0.78) {$\pi_i y$};
\draw [ <-] (0.1, 0.7) -- (0.7, 0.7);
\node at (-0.15, 1) {$H_i$};
\end{tikzpicture}
\caption{Projection $\pi_i y$ onto $H_i$ given by $\left\langle a_i, w\right\rangle = b_i$.}
\end{figure}
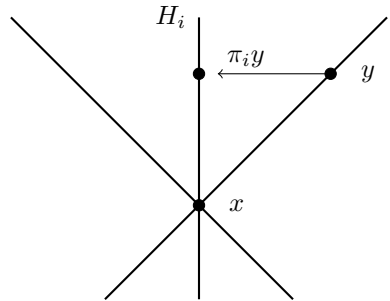
\end{center}

Using $\pi_i$ to denote the orthogonal projection $\pi_i: \mathbb{R}^n \rightarrow H_i$, one sees from Fig. 1 and the Pythagorean theorem  that 
$$ \forall w \in \mathbb{R}^n \quad~\forall 1 \leq i \leq m \qquad \|w - x\|^2 = \| \pi_i w - x \|^2 + \| w - \pi_i w\|^2 $$
from which we deduce
$$ \forall w \in \mathbb{R}^n \quad ~\forall 1 \leq i \leq m \qquad   \qquad \| \pi_i w - x \|^2 \leq  \|w - x\|^2.$$
This means that given an approximate solution $w$ of the linear system, projecting it onto any of the $m$ hyperplanes will always
lead to a better approximation of the solution. Moreover, each individual projection is cheap to compute since
$$\pi_i w = w+ \frac{b_i - \left\langle a_i, w\right\rangle}{\|a_i\|^2}a_i.$$
This has been first proposed as an iterative method for solving linear systems in the 1930s 
by Kaczmarz \cite{kac}, we refer to \S 2 for a more detailed discussion.

\subsection{Approximate Solutions}
We study approximate solutions $Ay \approx b$ which we define loosely as vectors $y \in \mathbb{R}^n$ for which
$\|Ay - b\|$ is small. We note that if $y$ is close to the true solution $x$, then 
$$ \| Ay - b\| = \| Ay - Ax \| \leq \|A\| \cdot \|x-y\|$$
and $y$ is an approximate solution. However, when $A$ is badly conditioned and has small singular values, then there are also approximate solutions that are far away from the true solution. Our main goal is to show that projection onto hyperplanes, a method that
has been proposed for finding solutions of the linear system $Ax = b$, is very good at finding approximate solutions. We start by stating the main result.

\begin{theorem}[Main Result] Let $A \in \mathbb{R}^{m \times n}$, $m \geq n$ have full rank, $\rk A = n$, and suppose $Ax = b$. Then, for any $x_0 \in \mathbb{R}^n$ and any $\varepsilon > 0$ there exists a sequence of $k$ projections of $x_0$ onto the hyperplanes $ x_k = \pi_{i_k} \pi_{i_{k-1}} \dots \pi_{i_1} x_0,$
where 
$$ k \leq 2 \log\left(\frac{\|A\|}{\varepsilon}\right) \cdot \frac{\|A\|^2_F}{\varepsilon^2}$$
such that $\|x_k - x\| \leq \|x_0 - x\|$ and
$$ \|A x_k - b\| \leq \varepsilon \cdot \| x_0 - x\|.$$
\end{theorem}

\textbf{Remarks.}
\begin{enumerate}
\item The upper bound on the number of projections $k$ only depends on the `size' of the matrix and the accuracy $\varepsilon$. It does not depend on the behavior of the small singular values and therefore does not depend on whether solving the linear system $Ax = b$ is difficult. This is expected: finding approximate solutions is easier than approximating solutions (by a factor proportional to the singular value in each subspace spanned by singular vectors).
\item Theorem 1 is an existence result. However, the proof suggests that if the projections are chosen randomly (with $\pi_i$ at each step being
chosen with likelihood proportional to $\|a_i\|^2/\|A\|_F^2$, this is the Strohmer-Vershynin Randomized Kaczmarz method \cite{strohmer}), then the arising sequence of projections satisfies the result with high likelihood (see \S 1.3 for a numerical example).
\item If $A$ is ill-conditioned, then there exist $x,y \in \mathbb{R}^n$ such that $Ax = b \approx Ay$ is a good approximate solution even though $\|x-y\|$ is very large: note that $\|x-y\| \leq \sigma_{\min}^{-1} \cdot \|Ax - Ay\|$ and badly conditioned matrices have very small singular values. This is a tremendous source of instability, approximate solutions $y$ could be many orders of magnitude larger than the true solution $x$. Since we have  $\|x_k - x\| \leq \|x_{k-1} - x\|$, the procedure is stable and
$$ \|x_k\| \leq 2 \|x\| + \|x_0\|.$$
\item The inequality is mainly of interest when $\varepsilon \geq \sigma_{\min}$. Once $\varepsilon < \sigma_{\min}$, one can use another argument (in \S 2.4) which shows that
$$ k \leq 2 \log\left(\frac{\|A\|}{\varepsilon}\right)  \frac{\|A\|_F^2}{\sigma_{\min}^2} \qquad \mbox{also suffices.}$$ 
 \end{enumerate}

Using the singular value decomposition, we see that there cannot too many singular values that are large: using $\sum_{i=1}^{n} \sigma_i^2 = \|A\|_F^2$, it follows that for all $\varepsilon > 0$
$$ \# \left\{ 1\leq i \leq n: \sigma_i \geq \varepsilon\right\} \leq \frac{\|A\|_F^2}{\varepsilon^2}$$
with equality if all singular values coincide. This means we can write
$ A = A_1 + A_2,$
where $\rk A_1 \leq \|A\|_F^2/\varepsilon^2$ has small rank and $\|A_2\| \leq \varepsilon$ has small operator norm.
In particular, the quantity $\|A\|_F^2/\varepsilon^2$ arises naturally as the effective numerical rank. Theorem 1, which is completely deterministic,
shows that it is possible to locate an approximate solution in a space of dimension $\sim \|A\|_F^2/\varepsilon^2$
using a number of projections that is both completely explicit (the rows of the matrix) and comparable 
(up to a logarithmic factor) to the numerical rank of the matrix.

\subsection{Numerical Example.}
We quickly illustrated Theorem 1 with an example.
The Hilbert matrix $H_n \in \mathbb{R}^{n \times n}$ is defined by
$$ (H_n)_{ij} = \frac{1}{i+j-1}.$$
It is symmetric, positive-definite and infamously ill-posed: its singular values decay geometrically and its condition number grows roughly like $\sim (1+\sqrt{2})^{4n}$. We set $n=1000$ and try to find approximate solutions of 
$$ H_{1000} x = b \qquad \mbox{where} \qquad b = H_{1000} (1,1,1,\dots, 1).$$
In particular, we know that $x = (1,1,\dots,1)$ but the algorithm does not.
Note that $H_{1000}$ is, albeit invertible, so very ill-posed that it is nearly impossible to determine solutions of $H_{1000}x = b$ unless we explicitly construct $b$ to have a simple solution $x$.

\begin{figure}[h!]
\begin{tikzpicture}[scale=1]
\node at (-0.2,0) {\includegraphics[width=0.43\textwidth]{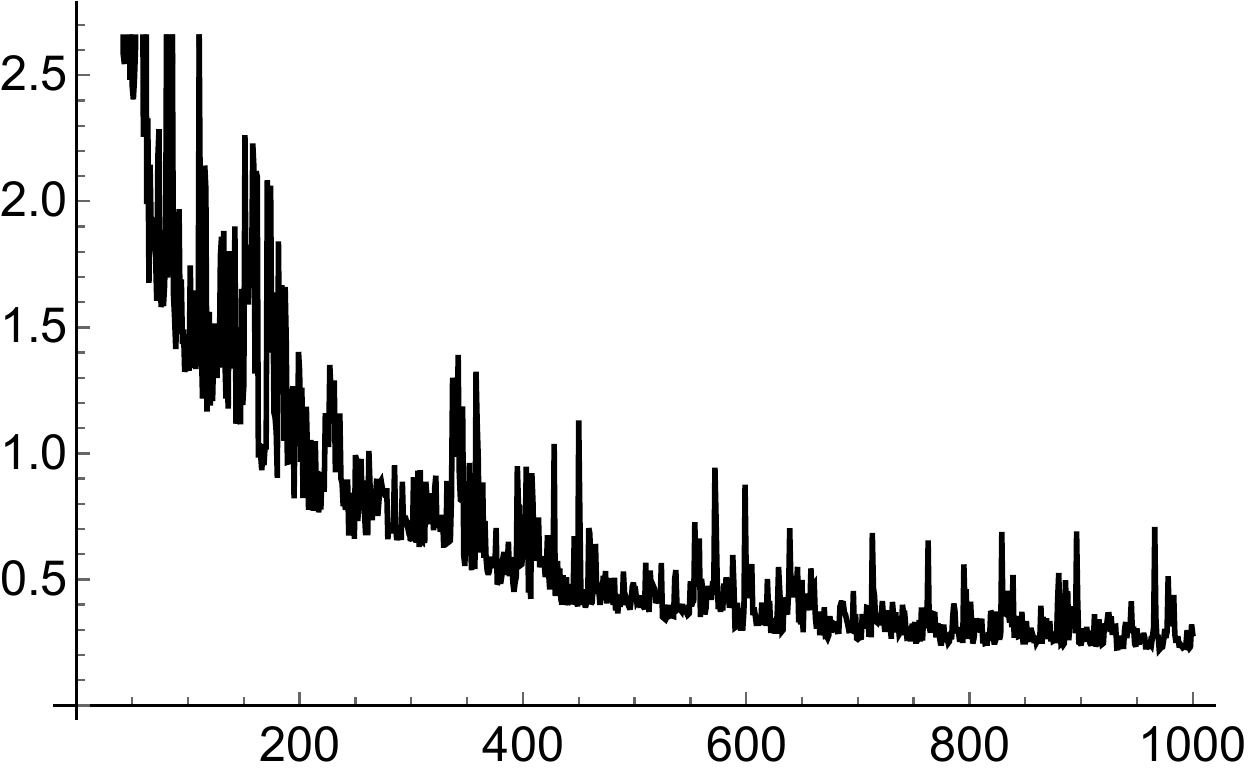}};
\node at (1, 1) {$ \|H_{1000}x_k - b\|$};
\node at (6-0.2,0) {\includegraphics[width=0.4\textwidth]{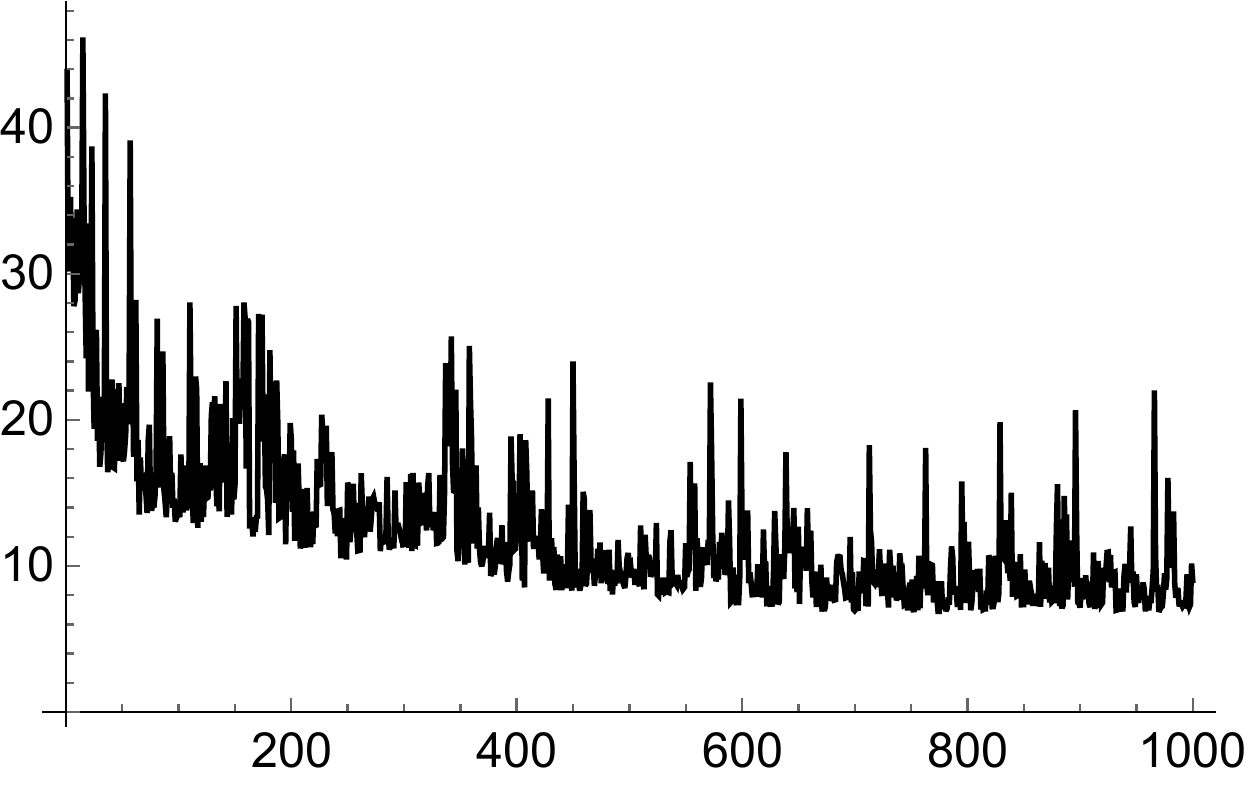}};
\node at (6+1, 1) {$\sqrt{k} \cdot \|H_{1000}x_k - b\|$};
\end{tikzpicture}
\captionsetup{width=0.9\textwidth}
\caption{Approximate solutions for the Hilbert matrix.}
\end{figure}
Our initial guess for a solution is $x_0 = (0,0, \dots, 0)$. We pick our sequence of projections randomly where $\pi_i$ is chosen with likelihood $\|(H_{1000})_i\|^2/\|H_{1000}\|_F^2$ in each step. Note that the initial error is $\| H_{1000}x_0 - b\| = \|b\| = \sqrt{1000} \sim 31$. We observe that $\|H_{1000} x_k -b\|$ decays dramatically (see Fig. 2). We also see that the decay seems to behave like $\sim 1/\sqrt{k}$ as suggested by our result.
The approximate solution $x_k$ can be seen in Fig. 3 as a function of its 1000 coordinates.  We see that $x_{10000}$ is not too far away from the ground truth (the constant vector $\mathbf{1}$). Note that $x_{10000}$ is a very good
approximate solution (see Fig. 3).
\begin{figure}[h!]
\begin{tikzpicture}
\node at (0,0) {\includegraphics[width=0.42\textwidth]{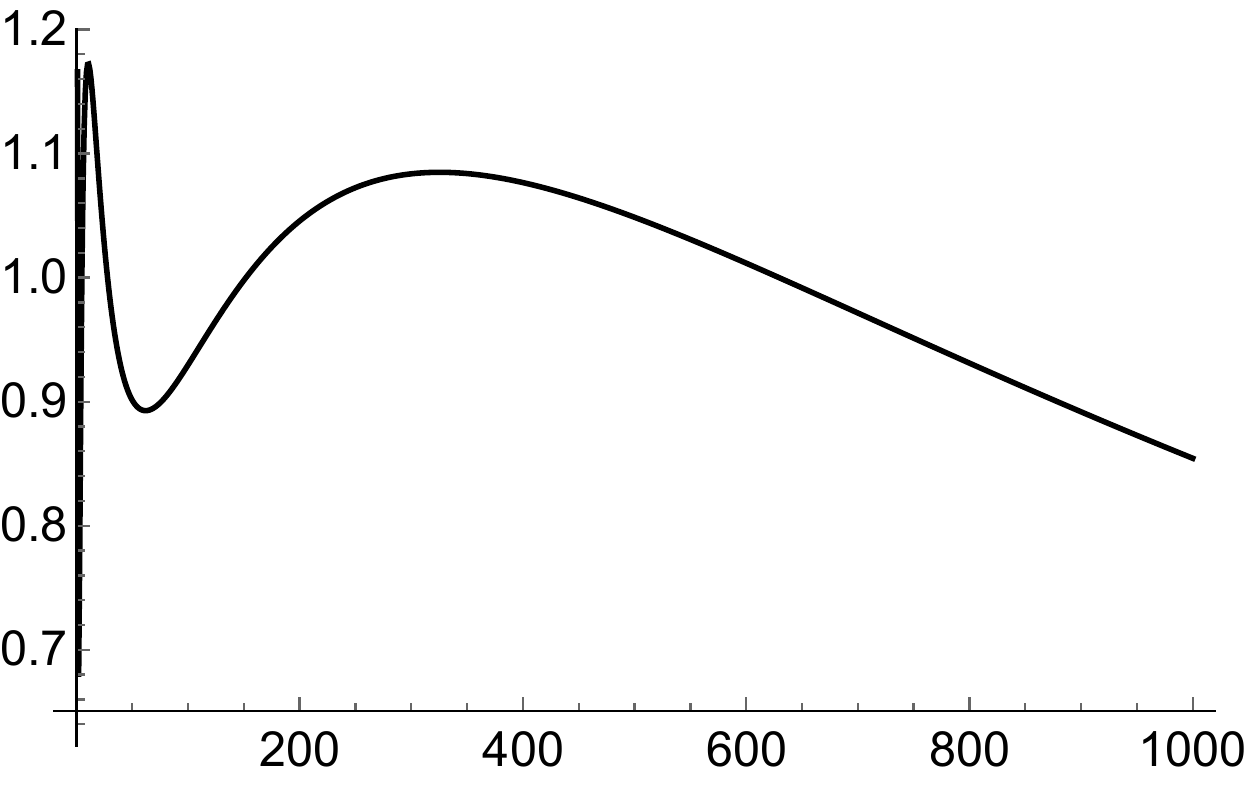}};
\node at (1, 1) {$x_{10000}$};
\node at (6.5,0) {\includegraphics[width=0.45\textwidth]{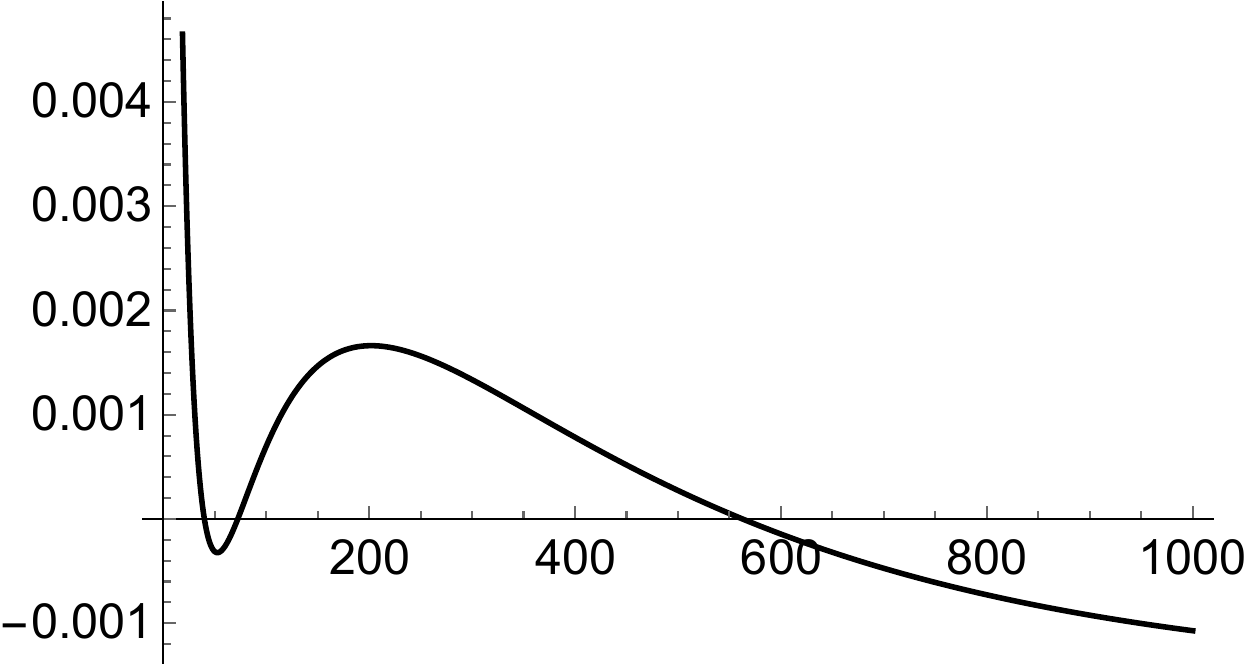}};
\node at (7.5, 0) {$  \|H_{1000} x_{10000} - b\|$};
\end{tikzpicture}
\captionsetup{width=0.9\textwidth}
\caption{$x_{10000}$ is a good approximate solution. The maximum error occurs in the first few entries and is small: $\|H_{1000} x_{10000} - b\|_{\infty} \sim 0.02$.}
\end{figure}

%
%
%
%

\section{The Geometry of Random Kaczmarz}
The main result of our paper is to show that there exists relatively short sequences of projections leading to good approximate solutions of linear system of equations.
Our proof makes uses of the Strohmer-Vershynin Random Kaczmarz method \cite{strohmer} which provides a framework in which the result can be proven while simultaneously providing a method for finding approximate solutions. We quickly summarize the method and its convergence properties (\S 2.1), how they
relate to Theorem 1 (\S 2.2) and discuss a new estimate for the method that will imply Theorem 1 (\S 2.3).

\subsection{The Kaczmarz method}  
 Let, as above, $A \in \mathbb{R}^{m \times n}$, $m \geq n$, be a (possibly overdetermined) linear system of equations, suppose that $A$ is injective and suppose that
$ Ax =b,$
where $x \in \mathbb{R}^n$ is the (unknown) unique solution and $b$ is a given right-hand side.
The Kaczmarz method proceeds as follows: given an approximation
of the solution $x_k$, take an arbitrary equation, say the $i-$th equation, and project $x_k$ onto the hyperplane described by
the $i-$th equation $ \left\langle a_i, x \right\rangle = b_i$, formally:
\begin{align}
 x_{k+1} = \pi_i x_k = x_k + \frac{b_i - \left\langle a_i, x_k\right\rangle}{\|a_i\|^2}a_i.
\end{align}
 The Pythagorean theorem implies that $\|x_{k+1} - x \| \leq \|x_k - x\|$ suggesting convergence of the method. It is not easy to make this quantitative \cite{gal}.
 There has been a substantial amount of work on this method, we refer to \cite{agmon, cenker, gal, kac, motzkin} as some of the fundamental papers in the area before 2009.  In 2009,
 Strohmer \& Vershynin \cite{strohmer} proposed to randomize the method and determined the convergence rate in terms of the smallest singular value and the Frobenius norm of $A$.

\begin{thm}[Strohmer \& Vershynin,  \cite{strohmer}] If the projection onto $H_i$ is chosen with likelihood $\|a_i\|^2/\|A\|_F^2$, then
\begin{align} \label{eq:strohmer}
\mathbb{E}~ \|x_k - x\|^2 \leq \left(1 - \frac{\sigma_{\min}^2}{\|A\|_F^2}\right)^k \|x_0 - x\|^2,
\end{align}
where $\sigma_{\min}$ is the smallest singular value of $A$.
\end{thm}

 The elegance of both the result and the proof has inspired a lot of subsequent work, see \cite{bai0, bai, gower, gower2, haddock, haddock15, haddock2, hadd, jiao, leventhal, ma, mar, need, need2, need25, need3, need4, nutini, stein2, stein3, stein4, stein5, zouz}.

\subsection{Approximate Solutions.} 
The expected projection of $x_k - x$ onto the smallest singular value decays exactly as predicted by Strohmer \& Vershynin while projections onto all
other singular vectors decay faster. 

\begin{thm}[Steinerberger, \cite{stein}] Let $v_{\ell}$ be a (right) singular vector of $A$ associated to the singular value $\sigma_{\ell}$. Then
$$\mathbb{E} \left\langle x_{k} - x, v_{\ell} \right\rangle = \left(1 - \frac{\sigma_{\ell}^2 }{\|A\|_F^2}\right)^k\left\langle x_0 - x, v_{\ell}\right\rangle.$$
\end{thm}

This has a couple of implications: in particular, singular vectors associated to the smallest singular value will asymptotically dominate. Moreover, even for moderately large $k$, we expect $x_k - x$ to be mainly a linear combination of singular vectors associated to small singular values. However, we also emphasize that this projection rather purely understood, we only know its value in expectation.
\begin{quote}
\textbf{Open Problem.} How does the random variable $\left\langle x_{k} - x, v_{\ell} \right\rangle$ behave?  How strongly is it concentrated around its mean?  Is it possible to get estimates on its higher moments?
\end{quote}

We will now proceed with a hypothetical computation assuming that there is some form of concentration around the mean and that 
\begin{align} \label{eq:conj}
\mathbb{E} \left\langle x_{k} - x, v_{\ell} \right\rangle^2 \quad \lesssim^{?}_{?} \quad  \left(1 - \frac{\sigma_{\ell}^2 }{\|A\|_F^2}\right)^{2k}\left\langle x_0 - x, v_{\ell}\right\rangle^2.
\end{align}
\eqref{eq:conj} is false: if $\left\langle x_0 - x, v_{\ell} \right\rangle = 0$, then this inequality would imply that $\left\langle x_k - x, v_{\ell} \right\rangle = 0$
for all $k \in \mathbb{N}$ which is clearly not the case. However, it is conceivable that something very close to \eqref{eq:conj} is true. We will use it as a convenient falsehood for the purpose of a napkin computation.
Assuming \eqref{eq:conj} and using the singular vectors as a basis
\begin{align*}
 \mathbb{E} \| A (x_k - x) \|^2 &= \mathbb{E}~ \sum_{i=1}^{n} \sigma_i^2 \left\langle x_{k} - x, v_{i} \right\rangle^2 \\
 &=   \sum_{i=1}^{n} \sigma_i^2  \cdot \mathbb{E} \left\langle x_{k} - x, v_{i} \right\rangle^2 \\
 &\lesssim^{?}_{\eqref{eq:conj}}  ~~ \sum_{i=1}^{n} \sigma_i^2  \left(1 - \frac{\sigma_{i}^2 }{\|A\|_F^2}\right)^{2k}\left\langle x_0 - x, v_{i}\right\rangle^2.
 \end{align*}
We use the elementrary inequality, valid for $0 < x < y$ and $z>0$,
$$ \left(1 -\frac{x}{y} \right)^z =  \left[ \left(1 -\frac{x}{y} \right)^{y}\right]^{\frac{z}{y}} \leq \exp\left( - \frac{ x z}{y} \right)$$
to deduce
$$  \mathbb{E} \| A (x_k - x) \|^2 \quad\lesssim^{?}_{\eqref{eq:conj}}  \quad \sum_{i=1}^{n} \sigma_i^2 \exp\left( - \frac{2 \sigma_{i}^2 k }{\|A\|_F^2} \right) \left\langle x_0 - x, v_{\ell}\right\rangle^2.$$
Combining this with the other elementary inequality, valid for $x, a > 0$, stating that
$x^2 \exp\left( - a x^2 \right) \leq  1/a$
we arrive at
$$  \mathbb{E} \| A (x_k - x) \|^2 \quad \lesssim^{?}_{\eqref{eq:conj}}  \quad \sum_{i=1}^{n} \frac{\|A\|^2_F}{k} \left\langle x_0 - x, v_{\ell}\right\rangle^2 \leq \frac{\|A\|_F^2}{k} \|x_0 - x\|^2.$$
If we want $\mathbb{E} \|Ax_k - b\| \leq \varepsilon \|x_0 - x\|$, then this would suggest we pick
$$  \frac{\|A\|_F}{\sqrt{k}} \leq \varepsilon \qquad \mbox{and thus} \qquad k \geq \frac{\|A\|_F^2}{\varepsilon^2}.$$
Thus the incorrect statement \eqref{eq:conj} would imply Theorem 1 without the logarithmic factor and underlines how a better understanding of the random variable $ \left\langle x_{k} - x, v_{\ell} \right\rangle$ might be useful. Needless to say, we will pursue a different route towards establishing Theorem 1.

\subsection{A Decay Estimate.} Our argument is based on a new kind of decay estimate for the standard Strohmer-Vershynin Random Kaczmarz method.

 \begin{theorem} Let $A \in \mathbb{R}^{m \times n}$, $m \geq n$ have full rank, $\rk A = n$, suppose $Ax = b$ and consider, for any given $x_0 \in \mathbb{R}^n$, a random sequence of vectors $(x_k)_{k=1}^{\infty}$ where $x_{k+1} = \pi_i x_k$ with likelihood $\|a_i\|^2/\|A\|_F^2$. Then, if no possible choice of $\leq k$ projections starting from $x_0$ lead to the exact solution, we have
 $$ \mathbb{E} \log \left( \frac{\|x_k-x\|}{\|x_0-x\|}  \right) \leq \frac{k}{2} \log\left(1 - \frac{\varepsilon^2}{\|A\|_F^2} \right) \cdot \frac{1}{k} \sum_{i=0}^{k-1}  \mathbb{P} \left( \frac{\|Ax_i -b\|}{\|x_i - x_0\|} \geq \varepsilon \right).$$
 \end{theorem}
 
A very natural choice is to take $\varepsilon$ to be the largest real number to ensure that the sum on the right-hand side evaluates to 1.  This number is easily seen to be, by definition, $\varepsilon = \sigma_{\min}$, the smallest singular value of the matrix $A$. Then $\|Ay\| \geq \sigma_{\min} \|y\|$ for all $y \in \mathbb{R}^n$ and the Theorem 2 can be rewritten as
 $$
   \mathbb{E} \log \left( \|x_k-x\|^2  \right) \leq k \log\left(1 - \frac{\sigma_{\min}^2}{\|A\|_F^2} \right) +  \log\left({\|x_0-x\|^2} \right)
 $$
which can be seen as a logarithmic version of the Strohmer-Vershynin bound \eqref{eq:strohmer}. Note that, for a positive random variable $X$, we have
$ \mathbb{E} \log(X) \leq \log \left( \mathbb{E} X \right)$ implying that the Strohmer-Vershynin bound is roughly one application of Jensen's inequality stronger. How big the difference between these estimates actually is depends on how strongly the random variable $X = \|x_k - x\|$ concentrates around its mean, another interesting problem to which we do not know the answer. For larger values of $\varepsilon > \sigma_{\min}$, we deduce that either decay happens more rapidly than expected or
that $\|Ax_i - b\|$ is a lot smaller than expected (either of which will be good for us).

\subsection{Comparison} One could argue that the Randomized Kaczmarz method by itself converges exponentially and should thus also lead to some sort of estimate. Using \eqref{eq:strohmer}, we have that iterates of the Randomized Kaczmarz method satisfy
$$  \mathbb{E}~ \| Ax_{k+1} - b\|^2 \leq \|A\|^2  \cdot \mathbb{E}~ \| x_{k+1} - x\|^2 \leq \|A\|^2 \left( 1 - \frac{\sigma_{\min}^2}{\|A\|_F^2} \right)^k \|x_0 - x\|^2.$$
We are looking for a vector satisfying
$ \|A x_k - b\| \leq \varepsilon \|x_0 - x\|,$
which is satisfied for
$$ k = 2 \log\left(\frac{\|A\|}{\varepsilon}\right)  \frac{\|A\|_F^2}{\sigma_{\min}^2}.$$
This estimate depends on the smallest singular value $\sigma_{\min}$ and thus does not have the same kind of stability as Theorem 1 which is independent of whether or not there are small singular values. This shows that a straight-forward application of \eqref{eq:strohmer} is not quite sufficient. It also shows that one can improve Theorem 1 in the regime $\varepsilon < \sigma_{\min}$ which, given Theorem 2, is not too surprising.

\begin{center}
\begin{figure}[h!]
\begin{tikzpicture}
\node at (0,0) {\includegraphics[width=0.5\textwidth]{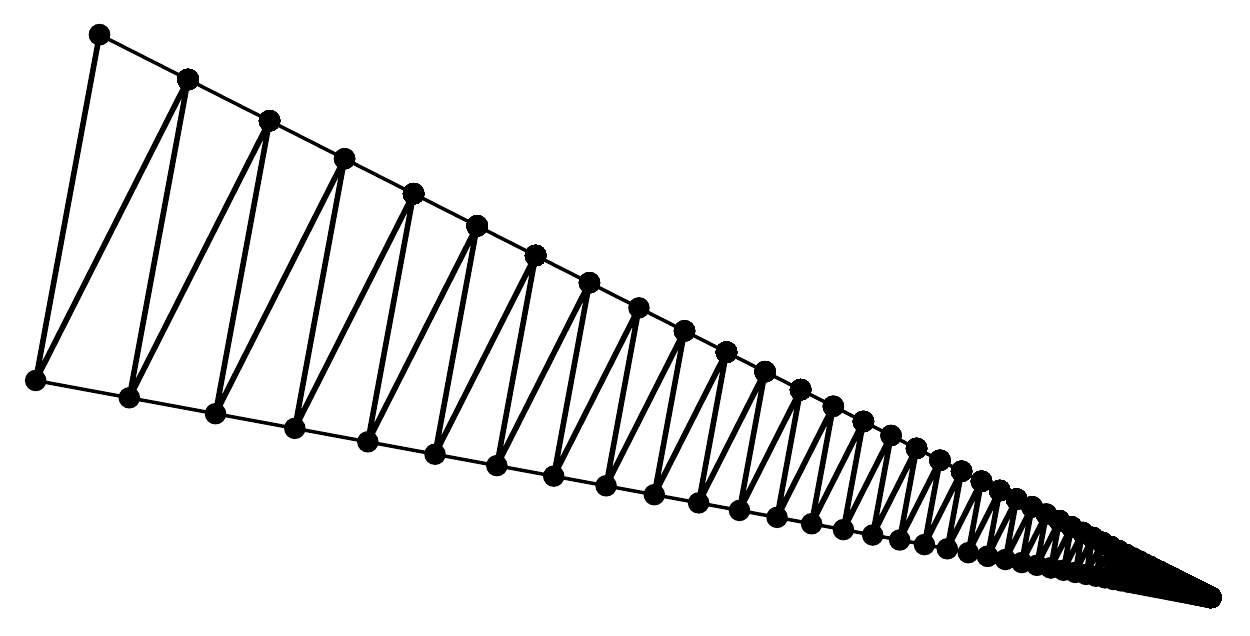}};
\node at (3.2, -1.6) {$x$};
\draw [thick] (4, 1) -- (8, -1);
\draw [thick] (4, 0.5) -- (8, -0.5);
\node at (8.3, -0.5) {$H_1$};
\node at (8.3, -1) {$H_2$};
\draw [thick,-> ] (6,0) -- (6.3, 0.75);
\draw [thick,-> ] (6,0) -- (5, 0.4);
\node at (7.5, 1) {\small leading singular vector};
\end{tikzpicture}
\caption{A cartoon sketch: slow convergence caused by hyperplanes meeting at a narrow angle. Such matrices have large singular vectors nearly orthogonal to the almost aligned lines while the second singular vector, pointing in the directions of the lines, corresponds to a small singular value: the matrix does not actually distinguish much between the bouncing approximations.}
\end{figure}
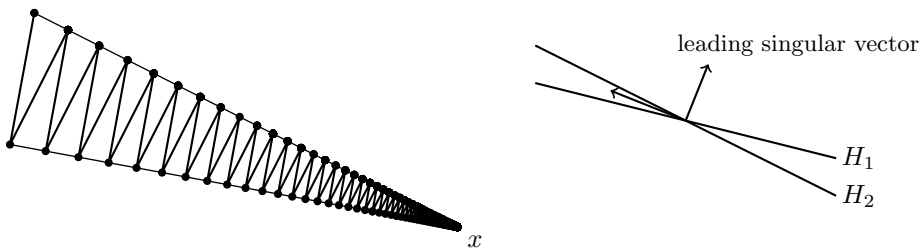
\end{center}

Another interpretation is as follows: the Random Kaczmarz method will exhibit slow convergence if there are small singular values.  Small singular values correspond to the hyperplanes being aligned in some fashion: this alignment simultaneously implies that the matrix will have large singular values acting orthogonal to the alignment and will not actually distinguish much between the aligned hyperplanes (meaning that there are singular vectors corresponding to small singular vectors that are approximately orthogonal to the normal vectors of the aligned hyperplanes). Theorem 1 and Theorem 2 imply that this happens in a scale-invariant fashion: all estimates are independent of the size of the small singular values.

\section{Proofs}

\subsection{Shrinking} We start with a fundamental geometric insight dating back to the seminal paper of Strohmer \& Vershynin \cite{strohmer}.
\begin{lemma} Let $x_k \in \mathbb{R}^n$ be arbitrary and choose $x_{k+1} = \pi_i x_k$
where the projection $\pi_i$ is chosen with likelihood proportional to $\|a_i\|^2/\|A\|_F^2$. Then
$$\mathbb{E} ~ \|x_{k+1} - x\|^2 = \left(1 - \frac{1}{\|A\|_F^2} \left\| A  \frac{x_k-x}{\|x_k - x\|} \right\|^2 \right) \|x_{k} - x\|^2.$$
\end{lemma}
The Lemma is not directly stated in Strohmer \& Vershynin \cite{strohmer} (but perhaps `stated in spirit'): arguing that $(x_k-x)/\|x_k - x\|$ is a vector of length 1, we note that
$$  \left\| A  \frac{x_k-x}{\|x_k - x\|} \right\|^2 \geq \sigma_{\min}^2  \left\|  \frac{x_k-x}{\|x_k - x\|} \right\|^2 = \sigma_{\min}^2$$
which then leads to the inequality commonly used in the literature. We are not aware of Lemma 1 being explicitly stated in the literature (though given how vast the literature has become and how simple the proof is, this is hard to rule out). A variant of Lemma 1 (phrased as an inequality) can be found in \cite[Theorem 2]{stein}. It will
be important for us that it actually is an equation.

\begin{proof} We have
\begin{align*}
\mathbb{E} ~ \|x_{k+1} - x \|^2 &= \frac{1}{\|A\|_F^2} \sum_{i=1}^{n} \|a_i\|^2 \left\|x_{k} - x -  \frac{\left\langle a_i, x_k \right\rangle}{\|a_i\|^2} a_i \right\|^2 \\
&= \|x_k -x\|^2 +  \sum_{i=1}^{n} \frac{ \|a_i\|^2 }{\|A\|_F^2} \frac{\left\langle a_i, x_k \right\rangle^2}{\|a_i\|^2} \\& - 2 \sum_{i=1}^{n} \frac{ \|a_i\|^2 }{\|A\|_F^2} \left\langle x_k -x, \frac{\left\langle a_i, x_k-x \right\rangle}{\|a_i\|^2} a_i  \right\rangle 
\\
&=  \|x_k-x\|^2 -  \sum_{i=1}^{n} \frac{1}{\|A\|_F^2} \left\langle a_i, x_k -x \right\rangle^2 \\
&= \|x_k-x\|^2 - \frac{1}{\|A\|_F^2}\sum_{i=1}^{n} \left\langle a_i, x_k-x \right\rangle^2 \\
&= \|x_k-x\|^2 - \frac{\|A(x_k -x)\|^2}{\|A\|_F^2} \\
&= \left(1 - \frac{1}{\|A\|_F^2}  \left\|  \frac{A(x_k-x)}{\|x_k-x\|} \right\|^2 \right)\|x_k-x\|^2.
\end{align*} 
\end{proof}

\subsection{Proof of Theorem 2} 
\begin{proof} Let $A$ and $x_0$ as well as $k$ be given. We may thus assume that they do not exist. In that case, we can write
$$ \|x_k  -x\| = \left( \prod_{i=0}^{k-1} \frac{\|x_{i+1} - x\|}{\|x_{i} -x\|} \right) \|x_0 - x\|$$
without the danger of dividing by 0. Taking a logarithm, we obtain
$$ \log (\|x_k-x\|) = \sum_{i=0}^{k-1}\log \left(  \frac{\|x_{i+1}-x\|}{\|x_{i}-x\|}  \right)  + \log(\|x_0-x\|).$$
The expectation of a sum of random variables is the sum of the expectations even if these random variables depend on each other, thus
$$ \mathbb{E} ~\log (\|x_k-x\|)  = \log(\|x_0-x\|) + \sum_{i=0}^{k-1} \mathbb{E} ~\log \left(  \frac{\|x_{i+1}-x\|}{\|x_{i}-x\|}  \right).$$
Appealing to Jensen's inequality, 
$$  \mathbb{E} ~\log \left(  \frac{\|x_{i+1}-x\|}{\|x_{i}-x\|}  \right) \leq   \log \left( \mathbb{E} \frac{\|x_{i+1}-x\|}{\|x_{i}-x\|}  \right)$$
which we use in conjunction with Lemma 1 to deduce
$$ \mathbb{E}~ \log \left(  \frac{\|x_{i+1}-x\|}{\|x_{i}-x\|}  \right) \leq  \mathbb{E} ~\frac{1}{2} \log\left(1 - \frac{1}{\|A\|_F^2} \left\| \frac{\|Ax_i - b\|}{\|x_i -x\|} \right\|^2 \right).$$
At this point, there is a case distinction: either $\|Ax_i - b\| \geq \varepsilon \|x_i - x\|$ (which happens with a certain likelihood) or not. This allows us to bound
$$ \mathbb{E}~ \log \left(  \frac{\|x_{i+1}-x\|}{\|x_{i}-x\|}  \right) \leq \frac{1}{2}  \log\left(1 - \frac{\varepsilon^2}{\|A\|_F^2} \right) \cdot \mathbb{P} \left( \frac{\|Ax_i - b\|}{\|x_i -x\|} \geq \varepsilon \right) $$
and collecting all the estimates leads to
$$ \mathbb{E} ~\log \left( \frac{\|x_k-x\|}{\|x_0-x\|}  \right) \leq \frac{k}{2} \log\left(1 - \frac{\varepsilon^2}{\|A\|_F^2} \right) \cdot \frac{1}{k} \sum_{i=0}^{k-1}  \mathbb{P} \left( \frac{\|Ax_i-b\|}{\|x_i-x\|} \geq \varepsilon \right) $$
which is the desired statement.
\end{proof}

\subsection{Proof of Theorem 1}
\begin{proof} Let $A, x_0$ and $\varepsilon > 0$ be given. We shall fix $k$ and then deduce a contradiction for $k$ sufficiently large.  If there is any possible choice of projections that would lead to the exact solution within $\leq k$ steps, Theorem 1 is trivially true. We can thus
suppose there does not exist selection of $\leq k$ consecutive projections such that $\|Ax_i -b\| \leq \varepsilon \|x - x_0\|$.  The Pythagorean Theorem implies the deterministic inequality
$ \|x_{i+1} - x\| \leq \|x_{i}- x\|$
and thus this implies there is also no possible selection of $\leq k$ projections such that
$$\|Ax_i -b\| \leq \varepsilon \|x - x_i\|.$$
Thus, independently of which projections are chosen, we have for all $i \leq k$ that
$$\mathbb{P} \left( \frac{\|Ax_i-b\|}{\|x_i-x\|} \geq \varepsilon \right) = 1.$$
Appealing to Theorem 2, we then have
\begin{align*}
 \mathbb{E} \log \left( \frac{\|x_k-x\|}{\|x_0-x\|}  \right) &\leq \frac{k}{2} \log\left(1 - \frac{\varepsilon^2}{\|A\|_F^2} \right) \cdot \frac{1}{k} \sum_{i=0}^{k-1}  \mathbb{P} \left( \frac{\|Ax_i-b\|}{\|x_i-x\|} \geq \varepsilon \right) \\
 &=  \frac{k}{2} \log\left(1 - \frac{\varepsilon^2}{\|A\|_F^2} \right).
 \end{align*}
A random variable satisfies $\mathbb{P}(X \leq \mathbb{E}X) > 0$.
Hence there exists at least one possible selection of $k$ projections such that the arising
$x_k = \pi_{i_k} \pi_{i_{k-1}} \dots \pi_1 x_0$, which we shall now fix, satisfies
$$ \log \left( \frac{\|x_k - x\|}{\|x_0 - x\|}  \right) \leq \frac{k}{2} \log\left(1 - \frac{\varepsilon^2}{\|A\|_F^2} \right).$$
This can be rewritten as 
$$ \frac{\|x_k - x\|}{\|x_0 -x\|} \leq \left(1 - \frac{\varepsilon^2}{\|A\|_F^2} \right)^{k/2}.$$
This implies that, since by assumption no selection of $\leq k$ hyperplanes leads to a good approximation, that
$$ \varepsilon \leq \frac{ \|A x_k  - b\|}{\|x_k - x_0\|} \leq \|A\| \frac{ \|x_k-x\|}{\|x_0-x\|}  \leq \|A\|  \left(1 - \frac{\varepsilon^2}{\|A\|_F^2} \right)^{k/2}.$$
This, in turn, leads to a contradiction as soon as
$$  \left(1 - \frac{\varepsilon^2}{\|A\|_F^2} \right)^{k/2} \leq \frac{\varepsilon}{\|A\|}.$$
Using the elementary inequality
$ \left(1-x/n\right)^n \leq e^{-x}$
and setting $k= 2 \ell \|A\|_F^2$, we arrive at a contradiction as soon as
$  \left(1 - \varepsilon^2/\|A\|_F^2 \right)^{k/2} \leq \exp\left(-\varepsilon^2 \ell\right).$
It now suffices to find an $\ell$ such that $\exp(-\varepsilon^2 \ell) \leq \varepsilon/\|A\|$ which requires
$$ \ell \geq \frac{1}{\varepsilon^2} \log\left(\frac{\|A\|}{\varepsilon}\right)$$
and thus we arrive at a contradiction for
$$ k = 2 \log\left(\frac{\|A\|}{\varepsilon}\right) \cdot \frac{\|A\|_F^2}{\varepsilon^2}.$$
\end{proof}

\end{document}